\documentclass[a4paper,11pt]{article}
\usepackage[english]{babel}
\usepackage[utf8]{inputenc}
\usepackage{paralist,url,verbatim, anysize}
\usepackage{amscd}
\usepackage{amsmath,amsthm,amssymb,amsfonts}
\usepackage{graphicx, graphics}
\usepackage[bookmarksopen=true]{hyperref}

\usepackage{color}
\theoremstyle{plain}
\newtheorem{theorem}{\bf Theorem}

\newtheorem{corollary}[theorem]{Corollary}
\newtheorem{lemma}[theorem]{Lemma}
\newtheorem{proposition}[theorem]{Proposition}
\newtheorem{mthm}{\bf Main Theorem}

\theoremstyle{definition}
\newtheorem{remark}[theorem]{Remark}

\newcommand{\sd}{\operatorname{sd} }

\newcommand{\R}{\mathbb{R}}

\newcommand{\CAT}{\mathrm{CAT}}

\begin{document}

\title{Linear embeddings of contractible and collapsible complexes}
\author{
Karim A. Adiprasito \\
\small Einstein Institute of Mathematics\\
\small Hebrew University, Jerusalem\\
\small \url{adiprasito@math.huji.ac.il}
\and
Bruno Benedetti\\
\small Department of Mathematics\\
\small University of Miami, Coral Gables\\
\small \url{bruno@math.miami.edu}
}

\date{August 4, 2020}
\maketitle 
\begin{abstract}
(1) We show that if a presentation of the trivial group is ``hard to trivialize'', in the sense that lots of Tietze moves are necessary to transform it into the trivial presentation, then the associated presentation complex (which is a contractible $2$-dimensional cell complex) is ``hard to embed in $\mathbb{R}^3$'', in the sense that lots of linear subdivisions are necessary. 

(2) For any $d$, we show that all collapsible $d$-complexes with $n$ facets linearly embed in $\R^{2d}$ after less than $n$ barycentric subdivisions. This is best possible, as cones over non-planar graphs do not topologically embed in $\R^{3}$.  
\end{abstract}

\section*{Introduction}

Every finite $d$-dimensional simplicial complex can be geometrically realized  in $\R^{2d+1}$ by placing its vertices in generic points.  
If we try to decrease  the ambient dimension by one, then not all $d$-complexes can be realized in $\R^{2d}$. For example, when $d=1$, not all graphs are planar. The class of planar graphs, completely characterized by Kuratowski's theorem, is invariant under subdivisions. 

In higher dimensions, the situation is much more complicated. When $d \ge 2$, realizing a complex or just a subdivision of it are no longer equivalent tasks. We say a complex 
\textit{linearly embeds} (resp.\ \textit{PL-embeds}) in $\R^k$, if the complex  (resp.\ some subdivision of it) can be realized in $\R^k$. We say that a complex \textit{topologically embeds} in $\R^k$, if it is homeomorphic to some subcomplex of $\R^k$.  Clearly, all linear embeddings are PL embeddings, and all PL embeddings are topological; but the converse inclusions are false.
By taking cones over non-planar graphs, one can easily produce examples of $d$-complexes that do not embed in $\R^{2d}$, neither linearly, nor PL, nor topologically. Whether a specific complex embeds in some $\R^k$ or not is a delicate matter; see e.g.~\cite{MTW} for a recent algorithmic approach. The conjecture that every contractible $2$-complex PL-embeds in $\R^{4}$ is wide open and related to the $4$-dimensional smooth Poincar\'e Conjecture \cite{Curtis}. (Weber proved that every acyclic $d$-complex PL-embeds in $\R^{2d}$ if $d \ge 3$ \cite{Weber}, but the argument does not extend to $d=2$.) In contrast, Horvati{\'c} and Kranjc showed that all acyclic $d$-complexes topologically embed in $\R^{2d}$ \cite{Hv, Kranjc}.

In the present paper, we prove the following result. Let us denote by $\Delta(a)$ a tower of exponentials of base $2$ of length $\lfloor a \rfloor$. In particular, $\Delta(0)=2$, $\Delta(1)=2^2$, and recursively $\Delta(n+1)=2^{\Delta(n)}$. We say a number is \emph{mild in $n$}, if it is bounded by a tower of exponentials in $n$ of bounded length.

\begin{mthm} \label{mthm:Contractible}
If a presentation of length $n$ of the trivial group needs at least $\Delta(\Omega(n))$ Tietze moves to trivialize, then the associated presentation complex does not embed linearly into $\R^3$ after a mild number of subdivisions. \end{mthm}

Here, $f=\Omega(g)$ is Knuth's notation, i.e., $f$ is asymptotically bounded from below by $g$. That presentations with the properties assumed in Main Theorem~\ref{mthm:Contractible} actually exist, is nontrivial; they were constructed by Bridson \cite{Bridson} and Lishak-Nabutovsky \cite{LN}. Main Theorem~\ref{mthm:Contractible} reflects the difference between dimension two and higher dimensions. In fact, in any dimension $d \ge 3$, the Freedman--Krushkal theorem claims that ``\emph{every topological embedding of a $d$-complex with $n$ facets in $\R^{2d}$ is isotopic to a linear embedding with $O(e^{n^{4+\varepsilon}})$ facets}'' \cite{FK}, a number which of course is much smaller than a tower of exponentials.

In the paper  \cite{FK} mentioned above, Freedman and Krushkal also show that for each $d \ge 2$ one can find $d$-complexes with $n$ facets that need a singly-exponential number of subdivisions to embed. We complement this result as follows:

\begin{mthm} \label{mthm:CollEmbed}
For any $d \in \mathbb{N}$, the $(n-1)$-st barycentric subdivision of any collapsible $d$-complex with $n$ facets
embeds linearly in $\R^{2d}$. 
\end{mthm}


\section*{Notation}
\paragraph*{Functions and their growth.} Given a function $f: \mathbb{N} \rightarrow \mathbb{N}$, we say that $f$ \emph{is mild} if it grows at most like a tower of exponentials of bounded length. For example, \[g(n)=2^{(3^{(4^n)})}\] is mild, because it grows like a tower of exponentials of length three. By $f=\Omega(g)$ we mean that $f$ is asymptotically bounded from below by $g$. By $\Delta(a)$ we denote a tower of exponentials of base $2$ of length $\lfloor a \rfloor$.

\paragraph*{Polytopal complexes.}
By $\R^d$ and $\mathbb{S}^d$ we denote the Euclidean $d$-space and the unit sphere in $\R^{d+1}$ with the standard (intrinsic) metric, respectively. A \emph{(Euclidean) polytope} in $\R^d$ is the convex hull of finitely many points in $\R^d$. A \emph{spherical polytope} in $\mathbb{S}^d$ is the convex hull of a finite number of points that all belong to some open hemisphere of $\mathbb{S}^d$. Spherical polytopes are in natural one-to-one correspondence with Euclidean polytopes, by taking radial projections. A \emph{geometric polytopal complex} in $\R^d$ (resp.\ in $\mathbb{S}^d$) is a finite collection of polytopes in $\R^d$ (resp.~$\mathbb{S}^d$) such that the intersection of any two polytopes is a face of both. For details, see e.g. the definition section in \cite{ABpart1}.

We denote by $\sd C$ the barycentric subdivision of $C$. Recursively, $sd^n (C) = \sd (\sd^{n-1} C)$. For the definitions of \emph{shellable}, \emph{line shelling}, etc. see e.g.~Ziegler \cite{Ziegler}. A \emph{free} face $\sigma$ is a face strictly contained in only one other face of $C$. An \emph{elementary collapse} is the deletion of a free face $\sigma$ from a polytopal complex~$C$.   We say that the complex $C$ \emph{collapses} to a subcomplex $D$, and write~$C\searrow D$, if $C$ can be reduced to $D$ by a sequence of elementary collapses. A \emph{collapsible} complex is a complex that collapses onto a single vertex. \emph{Non-evasiveness} can be defined by induction on the dimension, as follows: A simplicial complex of dimension $d$ is called \emph{non-evasive} if either $d=0$ and the complex consists of a single point, or if $d>0$ and the complex can be reduced to a single vertex by recursively deleting a vertex whose link is $(d-1)$-dimensional and non-evasive. 

\paragraph*{Group presentations.}
A presentation with generators $g_i$ and relators $r_j$ of a group $G$, usually denoted by
\[ G = < g_1, \ldots , g_s \: : \:\: r_1, \ldots, r_t > , \]
is called \emph{balanced} if $s=t$. 
Any balanced presentation $\wp$ can be associated with a cell complex $\mathfrak{C}=\mathfrak{C}(\wp)$, called \emph{presentation complex}, so that the $1$-cells of $\mathfrak{C}$ (resp. the $2$-cells of $\mathfrak{C}$) are in one-to-one correspondence with the generators (resp. the relators) of the presentation. If $G$ is the trivial group, then $C(\wp)$ is contractible; for details, see e.g.~\cite[Chapter I]{HMS}. The \emph{presentation length} $\ell(\wp)$ is the sum of all relator lengths in the presentation plus the number of generators. \emph{Tietze moves}  are the following four ways to change a presentation $\wp$ of a group $G$ into another presentation of the same group $G$:
\begin{compactenum}[I.]
\item Add a relation between generators that is a consequence of the existing ones.\\
For example, $ < x\, : \, x^3=1 >$  can be changed to $< x\, : \, x^3=1, x^6=1 >$.
\item Delete a relation between generators that is a consequence of the other ones.\\
For example, $< x\, : \, x^3=1, x^6=1 >$ can be changed to $ < x\, : \, x^3=1 >$.
\item Add a new generator that is expressed as a word in the other generators. \\
For example,  $ < x\, : \, x^3=1 >$ can be changed to $ < x, y \, : \, x^3=1, y=x^2>$.
\item If a relation can be formed where one of the generators is a word in the other generators, then remove that generator, replacing all occurrencies of it with the equivalent word. \\
For example, $ < x, y \, : \, xy=1, y=x^2>$ can be changed to $ < x\, : \, x^3=1 >$.
\end{compactenum}
The inverse of a type I move is a type II move; same for III and IV.

\bigskip

\section{Hard embeddings from hard presentations}
In this section we show that if $\wp$ is a complicated presentation of the trivial group, then some contractible $2$-complex $X(\wp)$ associated to it cannot embed in $\R^3$ in low complexity.
Our first step is to reduce ourselves to working with simplicial complexes:

\begin{lemma} \label{lem:1} Let $\wp$ be any finite presentation of a group $G$. Let  $\ell$ be the length of $\wp$. Then  the presentation complex $\mathfrak{C}(\wp)$ has
a triangulation $C=C(\wp)$ with at most $24\ell$ faces.
\end{lemma}

\begin{proof}
For each $i$, we realize the $2$-cell $\sigma_i$  corresponding to $r_i$ as a polygon with $\ell_i$ edges. This way we obtain a regular cell complex with $24\ell$ faces. Then we perform a barycentric subdivision.
\end{proof}

The second step is a revisitation of a classical construction due to Bing:

\begin{lemma}[{essentially Bing, cf.~\cite[Lemma 6]{Bing} and \cite[Theorem I.2A]{BingBook}}] \label{lem:2}
Let $X$ be any $k$-complex with $m$ faces 
that is geometrically realized in $\R^d$ (respectively, in $\mathbb{S}^d$). Then $X$ can be completed to a triangulation of a convex ball $B \subset \R^d$ (respectively, to a triangulation of $\mathbb{S}^d$) using  $O(m^{k})$ faces, for fixed $d$.
\end{lemma}

\begin{proof} We follow Bing's proof, and in particular, we use the terminology in \cite[Theorem I.2A]{BingBook}. We proceed by induction on $k$. When $k=1$ the claim is clear. Let $f_i(X)$ be the number of $i$-faces of $X$. 
We first ``shield off''  each $k$-face by taking the join with the boundary of a (suitably small) $(d-k)$-simplex. This introduces $f_k$ new $d$-dimensional simplices. Of course, $f_k \le m$. 
Once all $k$-faces are shielded off, we turn our attention to the link of each $(k-1)$-face. Such link is $k$-dimensional in $\mathbb{S}^{d-1}$, but all ``exposed faces" are of dimension $(k-1)$. By the inductive assumption, we can complete the triangulation of the link to a triangulation  of $\mathbb{S}^{d-1}$. This introduces $O(m^{k-1})$ faces for each $k$-face, for a total of $O(m^k)$ faces.
\end{proof}

\begin{corollary} \label{cor:2a}
Let $X$ be any $k$-complex with $m$ faces 
that is geometrically realized in $\R^d$ (or $\mathbb{S}^d$). There is a realization of the regular neighborhood $N_X$ of $X$ in $\R^d$ that
uses at most $O(m^{k})$ faces and collapses onto $N_X$.
\end{corollary}

\begin{proof}
As in Lemma~\ref{lem:2}, we complete $X$ to a triangulation $B$ of a convex ball, take two barycentric subdivisions, and look at the subcomplex of $B$ induced by all facets intersecting $X$. This gives the desired realization of the regular neighborhood. 
\end{proof}

Our third, crucial step shows how to extend a triangulation of a  $(d-1)$-sphere in $\R^d$ to a collapsible triangulation of its inside, using a mild number of faces.

\begin{proposition}  \label{lem:3}
Let $d \le 4$. Let $S$ be a triangulated $(d-1)$-sphere realized in $\R^d$ on $m$ faces. Let $B$ be the  ``inside'', i.e., the topological closure of the bounded connected component of the complement of $S$. There exist a triangulation of the $d$-ball $B$ such that
\begin{compactenum}[\rm (1)]
\item restricted to its boundary, $B$ coincides with a (mild) subdivision $S'$ of $S$;
\item $B$ has a mild number of faces;
\item $B$ is collapsible.
\end{compactenum}
\end{proposition}

\begin{proof} Every triangulated sphere of dimension $\le 2$ is polytopal. 
As for triangulated $3$-spheres, by a theorem of King \cite{King} polytopality can be achieved with a number of subdivisions estimated by $e^{O(m^2)}$. So without loss of generality we can assume that $S$ is polytopal. 
(Caveat: This does not mean that the inside of $S$ is convex.)

\begin{figure}[htb] \hskip-1mm
\includegraphics[width=0.43\linewidth]{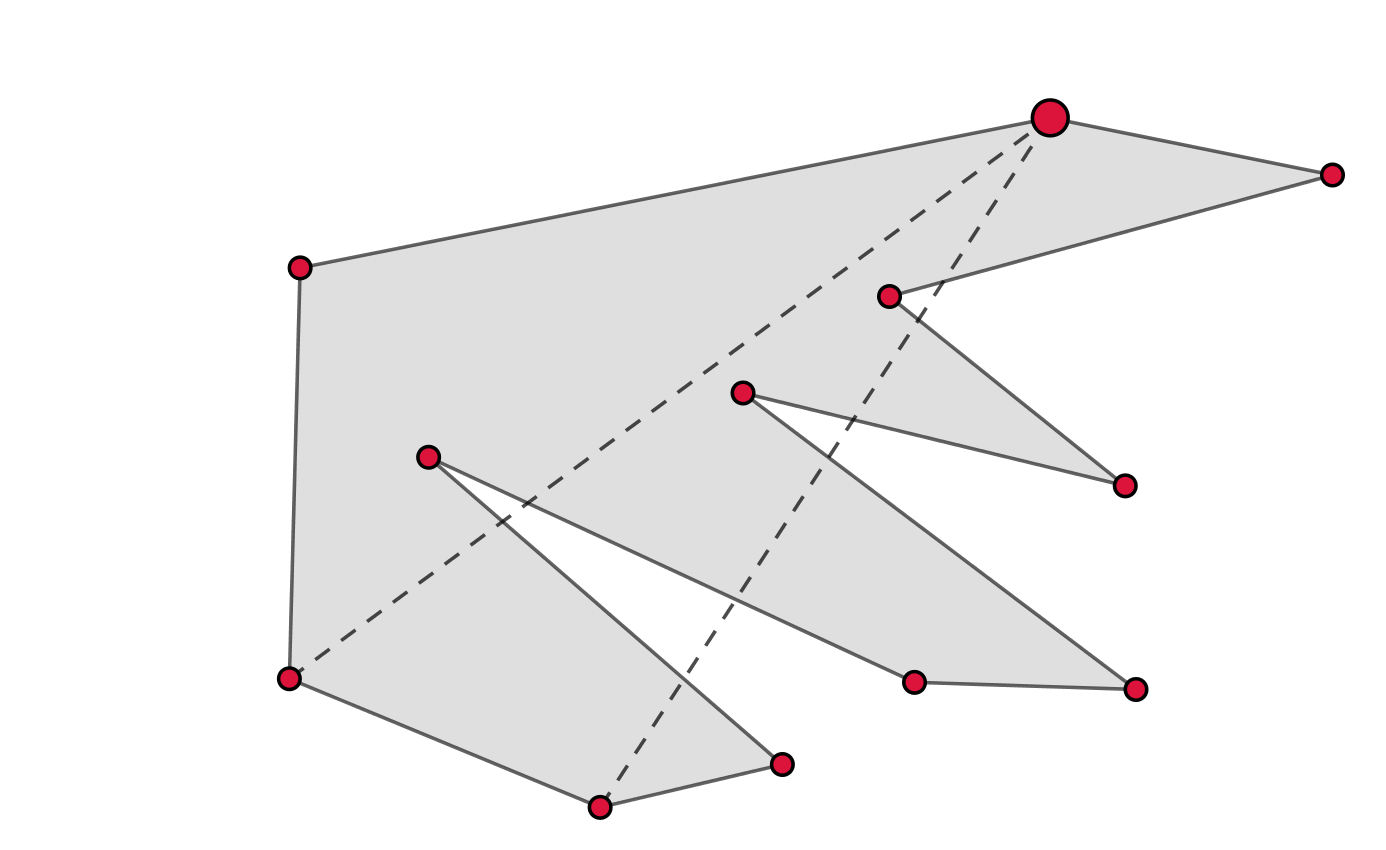} 
\includegraphics[width=0.43\linewidth]{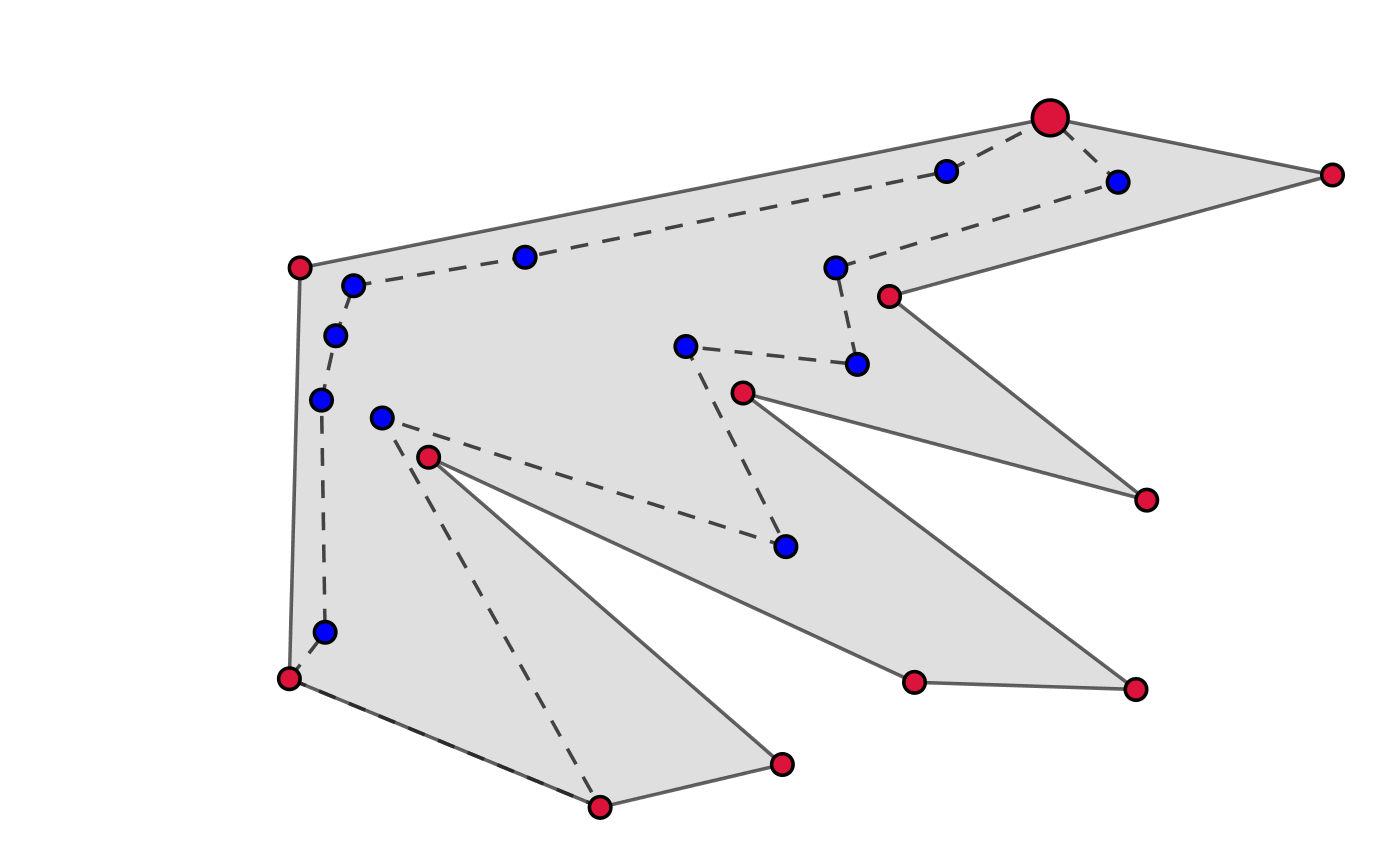} \hfill
\caption{Subdividing the geometric cone $C_1$ to ``make it fit'' inside $S$.}
\label{fig:1}
\end{figure}

Consider now a line shelling of $S$ that ends by removing the star of a vertex $v$. We want to triangulate the (abstract) cone over $S$ with apex $v$ in such a way that it forms a triangulation of $S$ in $\R^d$. For this, consider a first facet $F_1$ in the shelling order we selected. Let $C_1$ be the geometric cone over $F_1$ with apex $v$ in $\R^d$. 
Since the inside of $S$ is not necessarily convex, we do not know whether all facets of $S$ are visible from $v$. In particular, the geometric cone that we have just constructed might intersect $S$ in faces other than $F_1$ and $v$. To fix it, we want to move this cone to the inside of $B$. 

To achieve this, we may use a triangulation parallel to $S$ whenever the cone encounters the sphere $S$, using stellar subdivisions of $C_1$ at boundary faces (cf. Figure~\ref{fig:1} above.) The stellar subdivision of the cone is still collapsible, and we have introduced at most $e^{O(m^2)}$ new faces. 

\vskip1mm
Consider now the cone $C_2$ over the next facet $F_2$ with apex $v$. Again, we move $C_2$ to the interior by subdividing $(\partial C_2)\setminus C_1$ relative to $C_1$. Since the latter has $e^{O(m^2)}$ faces, this yields $e^{O(m^4)}$ new faces. 
This is repeated until the triangulation of the disk bounded by $S$ is completed, and gives us
a total of at most \[e^{m^{e^{O(m^2)}}}\] new faces. This number is mild in $m$, as desired.
\end{proof}

\begin{remark}
For polyhedra in $\mathbb{R}^3$ the bound above can be significantly improved, especially if one is willing to renounce to the collapsibility conclusion: See for instance \cite{ChaPa}.
\end{remark}

\begin{remark}
The proof of Proposition~\ref{lem:3} breaks down if $d\ge 5$: A triangulated $4$-sphere has a polytopal subdivision if and only if it is PL, and whether all $4$-spheres are PL is an important open problem, equivalent to the smooth Poincar\'e conjecture. But even if we restrict ourselves to PL spheres, making a PL $4$-sphere shellable may require a much larger number of barycentric subdivisions. In fact, let $F(n)$ be the smallest integer $k$ such that every PL $4$-sphere with $n$ facets becomes shellable after $k$ consecutive barycentric subdivisions. We know that $f(n)$ must be larger than a tower of exponentials, by the work of Lishak and Nabutovski \cite{LN, LNknots}. Were $f(n)$ bounded above by a computable function $h=h(n)$, then a shellability test for the $h(n)$-th barycentric subdivision would yield an algorithm to recognize PL $4$-spheres. It is conjectured that no such algorithm exists. If such conjecture holds true, then $F(n)$ must  be larger than any computable function of $n$.
\end{remark}

\begin{lemma}  \label{lem:4}
Let $D_1, D_2$ denote two geometric triangulations of the same $d$-manifold $D$ in $\R^d$, or in $\mathbb{S}^d$. Let $n_i$ be the number of faces of $D_i$ ($i=1,2$). Then the number of bistellar flips and stellar/inverse stellar subdivisions needed to connect $D_1$ and $D_2$ is mild in $n_1$ and in $n_2$.
\end{lemma}

\begin{proof}
Consider a regular subdivision $R$ of the convex hull of $D$ with no interior vertices. The number of faces of this triangulation is polynomial in the number $b$ of boundary vertices. In fact, by the Upper Bound Theorem \cite{Stanley}, it is bounded above by $O(b^{\left \lfloor \frac{d+1}{2} \right \rfloor})$. Let $r$ be the number of faces of $R$.

To find a regular stellar subdivision $D_1'$ of $D_1$ and $R$, we perform stellar subdivisions on $R$ for every transversal intersection of faces of $D_1$ and $R$, resulting in some regular subdivision $D_1''$ of $D_1$. This requires at most $n_1$ subdivisions. As a result, $D_1''$ has at most $r^{(d+1)^{n_1}}$ many faces, cf.~\cite[Lemma 3]{AIZ}.
This $D_1''$ is not necessarily stellar. However, any stellar subdivision of $D_1$ that refines $D_1''$ is necessarily regular: Compare the proof of \cite[Theorem 1]{AIZ}. This latter step requires at most $r^{(d+1)^{n_1}}$ subdivisions by \cite[Lemma 3]{AIZ}, and brings the count of faces to at most
\[n_1^{(d+1)^{r^{(d+1)^{n_1}}}}	\]
which is large but mild.
 Similarly, we find a common regular stellar subdivision $D_2'$ of $D_2$ and $R$ after a mild number of  subdivisions. The conclusion follows then from the fact that any two regular subdivisions of the same $d$-dimensional disk are connected by at most polynomially many bistellar flips in the number of faces, cf.~\cite[Corollary 5.3.11]{Triangulations}.
\end{proof}

\begin{lemma}[Lishak--Nabutovsky  \cite{LN}]  \label{lem:5}
Let $D_1$ and $D_2$ denote two geometric triangulations of the same convex subset $D$ in $\R^d$ or in $\mathbb{S}^d$. If $D_1$ and $D_2$ are connected by a single bistellar move, the induced presentations of their fundamental groups are at most $C_d$ Tietze moves apart, where $C_d$ depends only on $d$.
\end{lemma}

\begin{remark} \label{rem:6}
Any stellar subdivision can be achieved by polynomially many bistellar moves (in terms the number of the faces incident to the subdivided face). 
\end{remark}

\begin{proof}[\bf Proof of Main Theorem~\ref{mthm:Contractible}]
Consider a family of presentations of length $\ell$ that requires $\Delta(\Omega(\ell))$ Tietze moves to trivialize.
By Lemma \ref{lem:1}, in terms of $\ell$, the complex $C(\wp)$ has linearly many faces. Moreover, since $\wp$ presents the trivial group, the complex  $C(\wp)$  is contractible. 
Now, suppose there is a piecewise linear embedding $\varphi$; let us call $X$ the image of the complex $C(\wp)$ under $\varphi$. Assume by contradiction that the number of faces of $X$ is mild. Then there exists a mild triangulation of the regular neighborhood, which is a PL disk in $\R^3$ (being a contractible compact $3$-manifold). As observed, there exists a mild triangulation of that ball that is collapsible. Combining with the previous lemma, we obtain that the number of Tietze moves required to show triviality is mild; a contradiction.
\end{proof}

\begin{remark}
By Bridson \cite{Bridson} and Lishak-Nabutovsky \cite{LN}, there exists a family of presentations of the trivial group that require $\Delta(\Omega(\ell))$ Tietze moves to trivialize. So as direct corollary of Theorem~\ref{mthm:Contractible}, there are contractible 2-complexes $X$ with $n$ facets that do not linearly embed into $\mathbb{R}^3$ with a number of subdivisions  smaller than a tower of exponentials of length $n$.  This contrasts  a recent result by Freedman and Krushkal \cite{FK}, who proved a much smaller upper bound for the number of subdivisions necessary to PL embed complexes of dimension $d \ge 3$. 

Of course, if $Y$ is the cone over any non-planar graph, then $Y$ does not topologically embed in $\mathbb{R}^3$. The contractible complexes we build above though have a non-trivial, numerical obstruction to PL embeddability, in the sense that a priori, some might PL embed. 
\end{remark}

More generally, if $Y$ is the cone over any $(d-1)$-complex that does not topologically embed in $\mathbb{R}^{2d-2}$ (such as the $(d-1)$-skeleton of the $2d$-simplex), then $Y$ is a collapsible $d$-complex that does not topologically embed in $\mathbb{R}^{2d-1}$. In the next section, we complement this result by showing that if we try an ambient space of one dimension higher, then every collapsible complex linearly embeds after few subdivisions.

\section{All collapsible complexes embed after few subdivisions}
In this section we prove Main Theorem~\ref{mthm:CollEmbed}. 

\begin{lemma}\label{lem:facewise} Let $d, n$ be positive integers. Let $\Sigma$ be a $d$-simplex and let $\sigma$ be any of its facets.  Let $A_{n-1} = \sd^{n-1} (\partial \Sigma - \sigma)$. Then $\sd^{n} \Sigma$ admits a facewise linear map ${\varphi_n}$ to $A_{n-1}$  
that restricts to the identity on 
$\sd A_{n-1}$.
\end{lemma}

\begin{proof}
We claim that it is enough to prove this for $n=1$; in this case, $\varphi_1$ is defined by sending every face of $A_0=\partial\Sigma-\sigma$ to itself, and every remaining vertex in $\sd \Sigma$ to the unique vertex of $\Sigma$ not in $\sigma$. To recursively construct $\varphi_{n+1}$ when given $\varphi_{n}$, it suffices to observe that the barycentric subdivision of $A_{n-1}$ can be lifted to $\Sigma$ along the fibers of $\varphi_{n}$. See Figure~\ref{fig:2} below. 
\end{proof}

\begin{figure}[htb]
\centering
\includegraphics[width=0.8\linewidth]{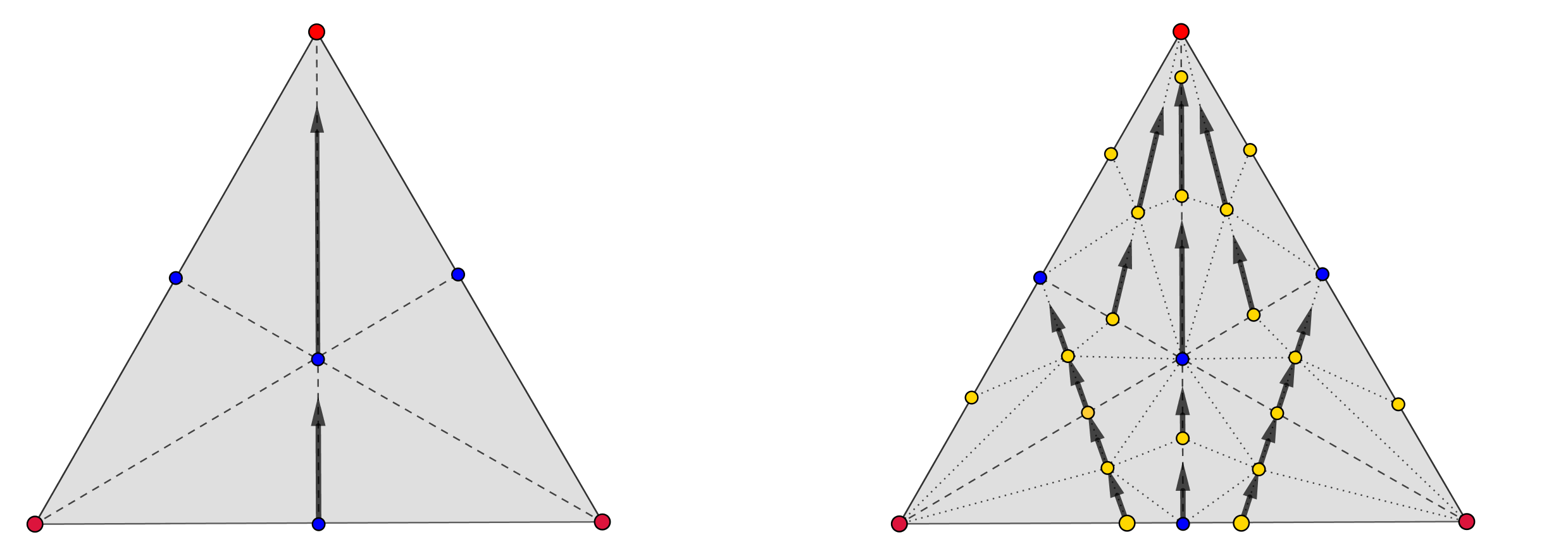}
\caption{How to construct $\varphi_0$ (left) and $\varphi_1$ (right).}
\label{fig:2}
\end{figure}

\begin{corollary}  \label{cor:compactness} 
 Let $d, n,k$ be positive integers. Let $\varepsilon>0$ be a real number. Let $\Sigma$ be a $d$-simplex and let $\sigma$ be one of its facets.  Let $A_{n-1} = \sd^{n-1} (\partial \Sigma - \sigma)$. Then any facewise linear embedding $\phi: A_{n-1} \longrightarrow \mathbb{R}^k$ can be extended to a facewise linear map $\phi:\sd^{n} \Sigma \rightarrow \mathbb{R}^k$ so that the image of $A_{n-1}$ and $\sd^{n} \Sigma$ differ by at most $\varepsilon$ in Hausdorff distance.
\end{corollary}

\begin{proof}
Using the map ${\varphi_n}$ of Lemma~\ref{lem:facewise}, we can deform $\Sigma$ until it gets close to $A_{n-1}$.
\end{proof}

\begin{proof}[\bf Proof of Main Theorem~\ref{mthm:CollEmbed}] We do not need to worry about the  lower-dimensional skeleton, since every $k$-complex generically embeds in $\mathbb{R}^{2d}$ when $k \le d-1$. Also, by genericity, we do not need to worry about intersections of $d$-faces with faces of dimension lower than $d$. Hence, we proceed by induction on the number $n$ of $d$-dimensional faces.

The case $n=1$ consists of a collapsible $d$-complex $C$ (not necessarily pure) with only one $d$-simplex. This certainly embeds in $\mathbb{R}^{2d}$ simply by generically embedding into $\mathbb{R}^{2d}$ the $(d-1)$-skeleton of $C$ and then by inserting the unique $d$-simplex. No barycentric subdivision is required.

Now, let $C$ be a collapsible $d$-dimensional simplicial complex with $n+1$ $d$-faces. Fix a collapse of $C$ and let $\sigma$ be the first free face removed in the sequence. Let $\Sigma$ be the unique $d$-face of $C$ containing $\sigma$. Set $C'=C - \sigma$. Since $C'$ is collapsible and has $n$ $d$-faces, by the inductive assumption we can find a geometric realization of $\sd^{n-1} C'$ in $\R^{2d}$. 
So all we need to do is to figure out how to position the vertices of the $n$-th barycentric subdivision of the simplex $\Sigma$ ``conveniently close'' to the geometric realization of $\sd^{n-1} C'$, so that self-intersections outside of stars of faces are avoided. This is precisely the task carried out by Corollary~\ref{cor:compactness}.
\end{proof}

\begin{remark} A $2$-dimensional complex $C$ is called \emph{$3$-thickenable} if there exists a triangulated $3$-manifold with boundary $M$ that collapses onto a subdivision of $C$. See Skopenkov \cite{Skopenkov} for the precise definition; we thank Uli Wagner for introducing us to this notion. Being $3$-thickenable is stronger than the planarity of all vertex links planar \cite{Skopenkov}. 
If $C$ is a $2$-dimensional complex that is $3$-thickenable, then reasoning exactly as in the proof of  Main Theorem~\ref{mthm:CollEmbed} above, one can conclude by induction on the number $n$ of facets that $C$  embeds linearly into $\mathbb{R}^3$.
\end{remark}

Applying the same reasoning to cubical complexes, we obtain an analogous result:

\begin{proposition} \label{prop:CAT0}
Every collapsible $d$-dimensional cubical complex with $n$ facets embeds linearly in $\R^{2d}$ after less than $n$ barycentric subdivisions.
\end{proposition}

\begin{corollary} \label{cor:CAT0}
Every $d$-dimensional (finite) $\CAT(0)$ cube complex with $n$ facets embeds linearly in $\R^{2d}$ after less than $n$ barycentric subdivisions.
\end{corollary}

\begin{proof}
By \cite[Corollary II]{ABpart1}, every $d$-dimensional  $\CAT(0)$ cube complex is collapsible.
\end{proof}

\section*{Acknowledgments}
We are most grateful to Martin Bridson, Boris Lishak, G\"unter Ziegler, Philip Brinkmann and Uli Wagner for useful conversations and for correcting earlier versions. 
Karim Adiprasito is supported by a Minerva fellowship of the Max Planck Society, an NSF Grant DMS 1128155, an ISF Grant 1050/16, and te grant ERC StG 716424 - CASe. Bruno Benedetti is supported by NSF Grants 1600741 and 1855165, a Provost Research Award, and a Scholarly \& Creative Activities Recognition Award by University of Miami. 
Part of this work was supported by the NSF under Grant No. DMS-1440140 while the authors were in residence at the Mathematical Sciences Research Institute in Berkeley, California (Fall 2017).

{\small

}

\end{document}